\documentclass[reqno,a4paper,11pt]{amsart}
\usepackage{a4wide}
\usepackage{amsthm,amsmath}
\usepackage[latin1]{inputenc}

\newcommand{\init}{\mathrm{in}}
\newcommand{\Kfm}{{K\left\langle X_1,\dots,X_n \right \rangle}}
\newcommand{\fm}{{\left\langle X_1,\dots,X_n \right \rangle}}
\newcommand{\G}{\mathrm{GL}(V)}
\newcommand{\gin}{\mathrm{gin}}
\def\setsuchas#1#2{\left\{\,{#1}\,\vrule\,{#2}\,\right\}}
\newcommand{\set}[1]{{\{#1\}}}
\newcommand{\divides}[2]{{#1 \left\lvert {#2} \right.}}
\newcommand{\dividesnot}[2]{{#1 \not \left\lvert {#2} \right.}}

\theoremstyle{definition}
\newtheorem{definition}{{Definition}}[section]
\newtheorem{example}[definition]{Example}

\theoremstyle{remark}
\newtheorem{note}[definition]{{Remark}}

\theoremstyle{plain}
\newtheorem{lemma}[definition]{{Lemma}}
\newtheorem{prop}[definition]{Proposition}

\newtheorem{theorem}[definition]{{Theorem}}

\newtheorem{corr}[definition]{{Corollary}}

\begin{document}

\title{Lifting Gröbner bases from the exterior algebra}
\author{Andreas Nilsson}
\address{Graduate School of Mathematical Sciences\\
                   University of Tokyo\\
    3-8-1 Komaba\\
 Meguro\\ Tokyo 153-8914\\
 Japan}
\email{andreas@ms.u-tokyo.ac.jp}

\author{Jan Snellman}
\address{Department of Mathematics\\
Link\"oping University\\
SE-58183 Link\"oping\\
SWEDEN
}
\email{jasne@mai.liu.se}


\sloppy
\begin{abstract}
  In \cite{Sturmfels:LiftGB}, the authors proved a number of results
  regarding Gröbner bases and initial ideals of those ideals \(J\),
  in the free associative algebra \(\Kfm\), which contain the
  commutator ideal. We prove similar results for ideals which contains 
  the anti-commutator ideal (the defining ideal of the exterior
  algebra). 

  We define one notion of generic initial ideals in \(\Kfm\), and show 
  that gin's of ideals containing the commutator ideal, or the
  anti-commutator ideal, are finitely generated.
\end{abstract}

\maketitle
\begin{section}{Notations}
  Let \(K\) be an infinite field, and let \(V\) be a vector space of
  dimension \(n\) over \(K\). Let \(X=\set{X_1,\dots,X_n}\) be an ordered
  basis of \(V\). Then the tensor algebra \(T(V)\) may be identified
  with the free associate (non-commutative) polynomial ring
  \(\Kfm\), which is the monoid
  ring on the free monoid 
  \(\left\langle X\right\rangle=\left\langle
    X_1,\dots,X_n\right\rangle\). 
  The symmetric algebra \(S(V)\) is then the quotient of \(\Kfm\) by
  the commutator ideal, the ideal generated by all \(X_i X_j - X_j
  X_i\) for \(i < j\). Alternatively, it is the monoid algebra on the
  free abelian monoid \([X]\).

  The exterior algebra over \(V\), \(E(V)\), can 
  be identified with the quotient of 
  \(K\left\langle X_1,\dots,X_n\right\rangle\) by the
  anti-commutator ideal \(C\) generated by \(X_iX_j + X_jX_i\), 
  \(1 \le i \le j \le n\). We denote the quotient epimorphism from
  \(\Kfm\) to \(E(V)\) by
  \(\pi\), and put  \(\pi(X_i) = x_i\). 
  If we denote the set of square-free products of elements in
  \(\pi(X)\) by \(Y\), then \(Y\) is a \(K\)-basis of \(E(V)\).
  For any subset \(I =\set{i_1,\dots,i_r} \subset \set{1,\dots,n}\), 
  \(i_1 < \cdots < i_r\),
  we denote by \(x_I\) the element \(x_{i_1} \cdots x_{i_r} \in
  Y\). If \(\sigma\) is any permutation of 
  \(\set{1,\dots,r}\) then \(x_{\sigma(i_1)} \dots x_{\sigma(i_r)} =
  \mathrm{sgn}(\sigma) x_I\). For \(j \in I\), we define
  \(\frac{x_I}{x_j}\) as \(x_{I \setminus \set{j}}\). Thus \(x_j
  \frac{x_I}{x_j} =  \genfrac{}{}{0pt}{3}{+}{-} x_I\).

  We let \(\prec_s\) be (the strict part of) a monoid
  well-order on \([X]\), and we let 
  \(\prec_e\) be its restriction to \(Y\). Then 
  \begin{equation}\label{eq:Yorder}
    m,n,t \in Y, \quad mt \neq 0, \, nt \neq 0, \quad m \prec_e n \qquad
    \implies \quad  mt \prec_e nt.
  \end{equation}

  \begin{definition}\label{def:order}
  We denote by \(<_{lex}\) the lexicographic order on \(\left \langle
    X \right \rangle\); it is
  supposed to order the variables 
  \[X_1 <_{lex} X_2 <_{lex} \cdots  <_{lex} X_n.\] 

  \noindent 
  We define a monoid well order on \(\left\langle X\right\rangle\) by
  \begin{equation}
    \label{eq:torder}
    M \prec_t N \quad \text{ if } 
    \begin{cases}
      0 \not \in \set{\pi(M), \pi(N)}, \quad \pi(M) \prec_e \pi(N)
      \\
      0 \not \in \set{\pi(M), \pi(N)}, \quad  \pi(M) = \pi(N), \quad 
     M <_{lex} N \\
      0 \in \set{\pi(M), \pi(N)},  \quad M <_{lex} N
    \end{cases}
  \end{equation}
  Hence, we see that 
\(X_iX_j \prec_t X_jX_i\) if \(i \le j\). 
  \end{definition}

  \begin{lemma}\label{lemma:CGB}
    The anti-commutators form a \(\prec_t\)-Gröbner basis for
    \(C\). Thus,
    \(\init_{\prec_t}(C) = \bigl(X_j X_i \, \lvert \, i \le j \bigr)\).
  \end{lemma}
  \begin{proof}
    It is enough to show that the anti-commutators form a
    \(<_{lex}\)-Gröbner basis. To do this, we observe that all
    obstructions resolve. 
  \end{proof}

  \begin{definition}\label{def:section}
We define a section of \(\pi\), 
  regarded as a \(K\)-linear map, as follows:
  \begin{equation}
    \label{eq:split}
    \begin{split}
      \delta: E(V) & \to T(V) \\
      x_{i_1} x_{i_2} \cdots x_{i_r} & \mapsto X_{i_1} X_{i_2} \cdots
      X_{i_r} \quad \text { if } \quad i_1 < i_2 < \cdots < i_r.
    \end{split}
  \end{equation}
    
  \end{definition}

  \begin{definition}\label{def:UL}
    Let \(L\) be a monomial ideal in \(E(V)\), and let
    \(L \ni m=x_{i_1}\cdots x_{i_r}\) with \(i_1 < \cdots <
    i_r\). Then we put 
    \begin{equation}
      \label{eq:UL}
      \mathcal{U}_L(m) = \setsuchas{u \text{  square-free monomial in 
          } x_{i_1+1}, \dots, x_{i_r -1}}{
        u \frac{m}{x_{i_1}} \not \in  L, \, 
        u \frac{m}{x_{i_r}} \not \in  L}
    \end{equation}
    If \(\mathcal{U}_L(m) = \set{1}\) for all generators \(m\) of
    \(L\), we say that \(L\) is \emph{squeezed}. 
  \end{definition}

  \begin{lemma}\label{lemma:Ufin}
    The set \(\mathcal{U}_L(m)\) is always finite. If \(\lvert i_r -
    i_1\rvert \le 1\) (in particular, if \(V\) is
    2-dimensional), then \(\mathcal{U}_L(m) = \set{1}\).
  \end{lemma}
  \begin{proof}
    There are only finitely many square-free monomials on \(x_{i_1+1},
    \dots, x_{i_r -1}\), and only the trivial such monomial if \(\lvert i_r -
    i_1\rvert \le 1\).
  \end{proof}

  We recall the definition of a \emph{stable ideal}: a monomial ideal \(L\)
  in \(E(V)\) is stable if, whenever \(m \in L\),
  \(\divides{x_j}{m}\), and \(\dividesnot{x_{j'}}{m}\) when \(j' > j\), 
  then it follows that \(x_i \frac{m}{x_j} \in L\) for all \(i \le j\).
  If the condition \(j'>j\) can be omitted, then \(L\) is
  \emph{strongly stable}. Clearly, strongly stable ideals are stable.
  \begin{lemma}\label{lemma:stabsq}
    Stable ideals are squeezed.
  \end{lemma}
  \begin{proof}
    Let \(L\) be a stable (square-free monomial) ideal, and let \(m\)
    be a generator, 
    \begin{displaymath}
      m= x_{i_1} \cdots x_{i_r}, \qquad 1 \le i_1 < \dots <
    i_r \le n.      
    \end{displaymath}
 Let furthermore 
 \begin{displaymath}
1 \neq u = x_{j_1} \dots x_{j_s}, \qquad i_1 < j_1 <
    \cdots < j_s < i_r. 
 \end{displaymath}
    Then the maximal index \(\ell\) such that
    \(\divides{x_\ell}{um}\) is \(\ell = i_r\). Hence stability implies that 
    \(x_{j_s}\frac{m}{x_{i_r}} \in L\), and since \(L\) is an ideal,
    \(\frac{u}{x_{j_s}}\cdot x_{j_s} \frac{m}{x_{i_r}} = 
    u \frac{m}{x_{i_r}} \in L\).  This shows that \(u \not \in
    \mathcal{U}_L(m)\).  
  \end{proof}
\end{section}

\begin{section}{The theorems}
  The following theorems are analogs of their counterparts in
  \cite{Sturmfels:LiftGB}.  
  For the remainder of this note, we use the following notation:
  \begin{center}
  \fbox{\(I \subset E(V)\) is a \emph{homogeneous
  ideal}, \emph{generated in degrees \(\ge 2\)}, and
\(J=\pi^{-1}(I)\). }
  \end{center}


  \begin{theorem}\label{thm:deltaGB}
    Let \(G\) be a minimal Gröbner basis of \(I\), with respect to
    \(\prec_e\). Then the anti-commutators 
    \(\setsuchas{X_iX_j +
      X_jX_i}{ 1 \le j \le i \le n}\), together with the set 
    \begin{equation}
      \label{liftGB}
      \setsuchas{\delta(u \cdot f)}{f \in G, \, u \in
        \mathcal{U}_{\init_{\prec_e}(I)}(\init_{\prec_e}(f))} 
    \end{equation}
    constitute a minimal Gröbner basis for \(J\), with respect to
    \(\prec_t\). 


    In particular, \(\init_{\prec_t}(J)\) is minimally generated by 
    \(\setsuchas{X_iX_j}{ 1 \le j \le i \le n}\) and
    \begin{displaymath}
      \label{eq:init}
      \setsuchas{\delta(u \cdot m)}{m \text{ minimal generator of }
        \init_{\prec_e}(I),\,  
        u \in \mathcal{U}_{\init_{\prec_e}(I)}(m)}.
    \end{displaymath}
  \end{theorem}
  \begin{proof}
    We start by noting that if
    \begin{displaymath}
      m_1 \in     \delta \left( \setsuchas{\init_{\prec_e}(g)}{g \in
          G} \right), \qquad
      m_2 \in \setsuchas{X_iX_j}{ 1 \le j \le i \le n},
    \end{displaymath}
    then \(m_1\) does not divide \(m_2\).
    This follows since
    \(G\) has no elements of degree 1, and since quadratic monomials
    are lifted to ``strictly ordered'' quadratic monomials.
    Hence \(\setsuchas{X_iX_j + X_jX_i}{ 1 \le j \le i \le n}\) is
    part of a minimal Gröbner basis, and \(\setsuchas{X_iX_j}{ 1 \le j
      \le i \le n}\) is part of a minimal generating system of
    \(\init_{\prec_t}(J)\). 

    Mimicking the proof in \cite{Sturmfels:LiftGB},
    we now show that
    if \(M = X_{i_1} \cdots X_{i_r} \in \fm\), then
    \(M \in \init_{\prec_t}(J)\) if and only if \(\pi(M) \in
    \init_{\prec_e}{I}\) or \(i_j \ge i_{j+1}\) for some \(j\).

    So, suppose first that  \(i_j \ge i_{j+1}\) for some \(j\). Then
    \(X_{i_j}X_{i_{j+1}} = \init_{\prec_t}(X_{i_j}X_{i_{j+1}} +
    X_{i_{j+1}}X_{i_j})\), and \(\divides{X_{i_j}X_{i_{j+1}}}{M}\),
    hence \(M \in \init_{\prec_t}(J)\). If on the other hand \(\pi(M) \in
    \init_{\prec_e}{I}\), but \(i_1 < \cdots < i_r\) (in this case, we 
    say that \(M\) is square-free and ordered), then
    \(M=\delta(x_{i_1} \cdots x_{i_r})\). Since \(\delta\) is a
    section to \(\pi\), \(\pi(M) = \pi(\delta(x_{i_1} \cdots x_{i_r})) =
    x_{i_1} \cdots x_{i_r}\). Hence, there exist an \(f \in I\) with
    \(\init_{\prec_e}(f) = x_{i_1} \cdots x_{i_r}\). Put \(F =
    \delta(f)\), then \(\pi(F) = f\), hence \(F \in J\). It is
    easy to see that \(\init_{\prec_t}(F) = M\). Therefore, \(M \in
    \init_{\prec_t}(J)\).

    For the converse, suppose that \(M \in \init_{\prec_t}(J)\). If
    there is some \(j\) such that \(i_j \ge i_{j+1}\), we are done;
    suppose therefore that \(i_1 < \cdots < i_r\). There exist some
    \(F \in J\) with \(\init_{\prec_t}(F) = M\). From the fact that
    \(i_1 < \cdots < i_r\) we conclude that \(\pi(M) \neq 0\), hence
    that \(\pi(F) \neq 0\). Therefore \(\pi(M) =
    \init_{\prec_e}(\pi(F)) \in \init_{\prec_e}{I}\).

    To move on with the proof, we observe that \(\init_{\prec_t}(J)\)
    contains all non-commutative monomials that fail to be both
    ordered and square-free,
    and that the initial monomials of the
    anti-commutators form a minimal generating set for these
    monomials.

    Now suppose that \(m' = u m\), where \(m = x_{i_1} \cdots
    x_{i_r}\) is a minimal generator of \(\init_{\prec_e}(I)\), with \(1 
    \le i_1 < \dots < i_r \le n\), and where \(u \in Y\) is any
    square-free monomial. It remains to show the following:
    \(\delta(um)\) is a minimal generator of \(\init_{\prec_t}(J)\) if 
    and only if \(u \in \mathcal{U}_{\init_{\prec_e}(I)}(m)\).

    For the ``only if'' part, suppose that \(\delta(um)\) is a minimal
    generator of \(\init_{\prec_t}(J)\). If \(u = 1\), then it belongs 
    to \(\mathcal{U}_{\init_{\prec_e}(I)}(m)\). Suppose therefore that
    there exists
    some \(j\) such that \(\divides{x_j}{u}\). We first show that
    \(j < i_r\). We can assume that \(j\) is the maximal index such
    that \(\divides{x_j}{u}\). 

    If it were the case that \(j \ge i_r\), then \(\delta(um) =
    \delta(m \frac{u}{x_j} x_j) = \delta(m 
    \frac{u}{x_j}) X_j\). Now \(\delta(m \frac{u}{x_j}) \in
    \init_{\prec_t}(J)\), which contradicts the fact that
    \(\delta(um)\) is a minimal generator of that ideal. 

    Secondly, we show that \(j > i_1\). We can here assume that \(j\)
    is the \emph{minimal}  index such
    that \(\divides{x_j}{u}\). If it were the case that \(j \le i_r\),
    then \(\delta(um) = 
    \delta(x_j \frac{u}{x_j} m) = X_j \delta(\frac{u}{x_j} m)\). 
    Now \(\delta(\frac{u}{x_j} m) \in
    \init_{\prec_t}(J)\), which contradicts the fact that
    \(\delta(um)\) is a minimal generator of that ideal. Hence \(j >
    i_1\). 

    So, we know that \(i_1 < j < i_r\), and thus that \(u\) is a
    square-free monomial in \(x_{i_1 
      +1}, \dots, x_{i_j -1}\). So \(\delta(um) = X_{i_1} \delta(u
    \frac{m}{x_{i_1}}) = \delta(u \frac{m}{x_{i_r}})X_{i_r}\), which
    implies that neither \(\delta(u \frac{m}{x_{i_1}})\) nor
    \(\delta(u \frac{m}{x_{i_r}})\) is in \(\init_{\prec_t}(J)\) (that 
    would contradict the fact that \(\delta(um)\) is a minimal
    generator), hence neither \(u \frac{m}{x_{i_1}}\) nor
    \(u \frac{m}{x_{i_r}}\) is in \(\init_{\prec_e}(I)\).

    For the ``if'' direction, suppose that \(u \in
    \mathcal{U}_{\init_{\prec_e}(I)}(m)\). That means that neither \(u
    \frac{m}{x_{i_1}}\) nor 
    \(u \frac{m}{x_{i_r}}\) is in \(\init_{\prec_e}(I)\), hence that
    neither \(\delta(u \frac{m}{x_{i_1}})\) nor 
    \(\delta(u \frac{m}{x_{i_r}})\) is in
    \(\init_{\prec_t}(J)\). Therefore, we can conclude that
    \(\delta(um)\) is a minimal generator of \(\init_{\prec_t}(J)\).

  \end{proof}

  \begin{note}
    If \(I\) has linear generators, the set considered above isn't
    necessary a minimal Gröbner basis. As an 
    example, let \(V\) be 2-dimensional, 
    and let \(I\) be the principal ideal on
    \(x_1\). Then \(\pi^{-1}(I)=J=(X_1,X_1^2,X_1X_2+X_2X_1,X_2^2)=
    (X_1,X_2^2)\), so two of the anti-commutators are  redundant as
    generators. The same example shows that in Theorem 
    2.1 of \cite{Sturmfels:LiftGB}, we must disallow linear generators 
    of the ideal. 
  \end{note}

  \begin{corr}\label{corr:sq}
    A minimal Gröbner basis of \(I\) lifts to a minimal Gröbner basis
    of \(J\) (in other words, the lift of its elements, together with
    the anti-commutators, constitute a minimal Gröbner basis) if and
    only if \(\init_{\prec_e}(I)\) is squeezed.
  \end{corr}

  \begin{corr}\label{corr:fg}
    \(J\) has a finite Gröbner basis.
  \end{corr}
  \begin{proof}
    Follows from Lemma~\ref{lemma:Ufin} and Theorem~\ref{thm:deltaGB}. We
    can also see this directly: since the exterior algebra has finite
    vector space dimension, so has \(T(V)/J\), hence so has
    \(T(V)/\init_{\prec_t}(J)\). Hence, for some \(d\),
    \(\init_{\prec_t}(J)\) contains all monomials of total degree
    \(\ge d\). Hence, \(\init_{\prec_t}(J)\) is generated in degrees
    \(\le d\) and is therefore finitely generated.
  \end{proof}

  \begin{theorem}\label{thm:ginlift}
    After a generic linear change of coordinates \(g\), \(g(J)\) has a 
    Gröbner basis which is the ``lift'' of a Gröbner basis of \(g(I)\).
  \end{theorem}
  \begin{proof}
    It is proved in \cite{Aramova:Gotzman} that (regardless of the
    characteristic of \(K\)) the generic initial ideal of \(I\) is
    strongly stable. Hence, for a generic \(g \in \mathrm{GL(V)}\), 
    \(\init_{\prec_e}(g(I))\) is strongly stable, hence stable, hence
    squeezed. By Corollary~\ref{corr:sq}, \(\pi^{-1}(g(I))\) has a
    Gröbner basis which is the ``lift'' of a Gröbner basis of
    \(g(I)\). Since the anti-commutator ideal is invariant
    under linear changes of variables, we have that \(\pi^{-1}(g(I)) = 
    g(\pi^{-1}(I)) = g(J)\).
  \end{proof}

  \begin{note}
    This is the only time that we need the term-order \(\prec_e\) to
    be the restriction of a term order on \([X]\) (in the terminology
    of \cite{BoolTerm}, \(\prec_e\) is a \emph{coherent} boolean term
    order). The results in \cite{Aramova:Gotzman} only cover such term 
    orders.

    For the rest of our results, it is enough that \(\prec_e\) is a
    total order satisfying \eqref{eq:Yorder}.
  \end{note}
\end{section}

\begin{section}{Generic initial ideals in the tensor algebra}
  We denote the set of invertible linear transformations on \(V\) by
  \(\G\). Since \(\G\) acts (on the left) on \(V\), it also acts on
  the tensor algebra \(T(V)\). When \(T(V)\) is identified with
  \(\Kfm\), and \(G\) with the set of invertible \(n\) by \(n\)
  matrices, the action of \(g=\left(g_{ij}\right)_{i,j=1}^n\)
  on the non-commutative monomial \(X_{i_1} \cdots X_{i_r}\) is given
  by
  \begin{equation}
    \label{eq:action}
    g(X_{i_1} \cdots X_{i_r}) = g(X_{i_1}) \cdots g(X_{i_r}) =
    \left( \sum_{\ell = 1}^n g_{\ell, i_1} X_\ell \right) \cdots
    \left( \sum_{\ell = 1}^n g_{\ell, i_r} X_\ell \right).
  \end{equation}
  Note that since the commutator ideal, and the anti-commutator ideal, 
  are invariant under this action, there is an induced action of
  \(\G\) on the symmetric algebra and on the exterior algebra. The
  theory for generic initial ideals in these algebras, with respect to 
  these actions, is well known \cite{Aramova:Gotzman, GreStill:GIN,
    Green:gin}, but we belive that the construction we present here is 
  new. 

  We give \(\G\) the Zariski topology. An \(F_\delta\)-subset of
  \(\G\) is then a countable intersection of open sets. 
   {\underline{In what follows, we assume
    that \(K=\mathbb{C}\)}};  then such subsets are non-empty and
  dense, so a   condition which holds on an \(F_\delta\) subset holds
  for a   ``randomly choosen'' element.

  \begin{theorem}\label{thm:ginextass}
    Let \(\mathfrak{a} \subset \Kfm\) be a two-sided homogeneous
    ideal, and let 
    \(>\) be a term order on the free abelian monoid
    \(\fm\). Then there is a monoid ideal \(\mathfrak{b} \subset
    \Kfm\) and a \(F_\delta\)-subset \(U \subset \G\) such that 
    \begin{equation}
      \label{eq:noncomU}
      U \ni g \implies \init_>(g(\mathfrak{a})) = \mathfrak{b}
    \end{equation}
    If \(\mathfrak{b}\) is finitely generated then \(U\) is open.
  \end{theorem}
  \begin{proof}
    The proof is almost word-for-word identical with the proof in
    \cite{Ebud:View} of the corresponding assertion for the
    symmetrical algebra. To stress the similarity, we temporarily
    denote \(\Kfm\) by \(S\). We start by fixing a total degree \(d\)
    and a basis \(f_1,\dots,f_t\) for the \(K\)-vector
    space \(\mathfrak{a}_d\). Then, we let \(f=f_1 \wedge \cdots \wedge f_t
    \in \bigwedge^t S_d\), noting that there is a 1-1
    correspondence between \(t\)-dimensional subspaces of \(S_d\) and
    \(1\)-dimensional subspaces of \(\bigwedge^t S_d\). Furthermore,
    every element in \(\bigwedge^t S_d\) can be written as a finite
    sum of terms \(c m_1 \wedge \cdots \wedge m_t\), where \(c \in K\)
    and \(m_i\) are (non-commutative) monomials of total degree \(d\).
    We can assume that such a term is written such that \(m_1 > \cdots
    > m_t\). Then, we order such expressions ``lexicographically'',
    that is,
    \begin{displaymath}
      c m_1 \wedge \cdots \wedge m_t > d m_1' \wedge \cdots \wedge m_t'
    \end{displaymath}
    iff there is an \(i\) such that \(m_j = m_j'\) for \(j < i\), and
    \(m_i > m_i'\) (here of course we use the term order \(>\) for
    comparison). 
    
    With this convention, we can define the initial term \(\init(f)\)
    of \(f\). It is easy to see that in fact 
    \begin{displaymath}
    \init(f) = 
    \begin{cases}
      \init(f_1) \wedge \cdots \wedge \init(f_t) & \text{ if all
        \(\init(f_i)\) are different,} \\
      0 & \text{ otherwise.}
    \end{cases}
  \end{displaymath}
  We may perform ``Gaussian elimination'' on the \(f_i\), so that
  their linear hull remains unchanged, but their initial terms become
  distinct. 

    Now let \(h\) be an \(n \times
    n\) matrix of indeterminates, \(h = (h_{ij})\), and let \(h\) act
    on \(f\) by \(h(f) = h(f_1) \wedge \cdots \wedge h(f_t)\), where
    on each component it acts as  an element of \(\G\).
    
    Now, let \(m = m_1 \wedge \cdots \wedge m_t\) be the largest (with
    respect to the order just defined) term of \(h(f)\) with non-zero
    coefficient \(p=p(h_{ij})\). Define \(U_d \subset \G\) to be the
    open subset consisting of all \(g =(g_{ij}) \in \G\) such that
    \(p(g_{ij}) \neq 0\). Define \(\mathfrak{b}_d\) to be the \(K\)-vector
    subspace of \(S_d\) generated by \(m_1,\dots,m_t\).  Then
    \(\init(g(\mathfrak{a}))_d = \mathfrak{b}_d\) iff \(g \in U_d\).
    
    Define \(\mathfrak{b} = \bigoplus_d \mathfrak{b}_d\). To see that
    \(\mathfrak{b}\) is a two-sided ideal, it 
    suffices to show that \(S_1 \mathfrak{b}_d \subset
    \mathfrak{b}_{d+1} \supset \mathfrak{b}_d S_1\). 
    Since both \(U_d\) and \(U_{d+1}\) are open and non-empty, their
    intersection is also open and non-empty. Pick a \(g \subset U_d
    \cap U_{d+1}\). Then \(\mathfrak{b}_d = \init(g(\mathfrak{a}))_d\)
    and \(\mathfrak{b}_{d+1} = 
    \init(g(\mathfrak{a}))_{d+1}\). Since \(\init(g(\mathfrak{a}))\)
    is a two-sided ideal, we have 
    that 
    \[S_1 \init(g(\mathfrak{a}))_d \subset
    \init(g(\mathfrak{a}))_{d+1} \supset 
    \init(g(\mathfrak{a}))_d S_1,\]
    hence that
    \(S_1 \mathfrak{b}_d \subset \mathfrak{b}_{d+1} \supset
    \mathfrak{b}_d S_1\). 
    
    Put \(U = \bigcap_d U_d\).  Clearly, \(\init(g(\mathfrak{a})) =
    \mathfrak{b}\) iff \(g 
    \in U\). If \(\mathfrak{b}\) happens to be generated in degrees \(\le e\),
    then we claim that \(U = \bigcap_{d=1}^e U_d\), hence \(U\) is open.
    Suppose that \(g \in \bigcap_{d=1}^e U_d\), then
    \(\init(g(\mathfrak{a}))_d = 
    \mathfrak{b}_d\) for all \(d \le e\). Since \(\mathfrak{b}\) was
    supposed to be 
    generated in degrees \(\le e\), this implies that \(\init(g(\mathfrak{a}))
    \supset \mathfrak{b}\). But for all \(d\), \(\mathfrak{b}_d\) and
    \(\init(g(\mathfrak{a}))_d\) have 
    the same dimension as \(K\)-vector spaces, namely \(\dim_K
    \mathfrak{a}_d\). 
    Hence we conclude that in fact \(\init(g(\mathfrak{a})) = \mathfrak{b}\).
  \end{proof}

  We denote the ideal \(\mathfrak{b}\) above by \(\gin(\mathfrak{a})\)
  and call it the 
  \emph{generic initial ideal} of \(\mathfrak{a}\).
  \begin{lemma}\label{lemma:H}
    \(\gin(\mathfrak{a})\) and \(\mathfrak{a}\) have the same Hilbert series.
  \end{lemma}
  \begin{proof}
    Hilbert series are preserved by non-singular linear changes of
    coordinates, and by passing to the initial ideal.
  \end{proof}

  \begin{prop}
    \(\gin(\mathfrak{a})\) need not be finitely generated, even if
    \(\mathfrak{a}\) is.
  \end{prop}
  \begin{proof}
    There exists a finitely presented algebra with non-rational
    Hilbert series \cite{Shearer:NonRat}. By Lemma~\ref{lemma:H} the
    gin of the defining ideal must have non-rational Hilbert series,
    too. Since finitely generated monomial ideals have rational
    Hilbert series \cite{Ufnarovski:Criterion}, the result follows.
  \end{proof}

  However, combining the results of \cite{Sturmfels:LiftGB} with
  Theorem~\ref{thm:ginlift} we get that 
  \begin{theorem}\label{thm:fingin}
    If \(\mathfrak{a}\) contains the commutator ideal, or the
    anti-commutator ideal, but no non-zero linear form, then
    \(\gin(\mathfrak{a})\) is minimally generated
    by \(\delta\) of the  
    minimal generators of \(\gin(\pi(\mathfrak{a}))\), together with
    \(X_iX_j\), \(i > j\) (\(i \ge j\) for the exterior algebra).
  \end{theorem}
  Here, \(\pi\) and \(\delta\) are either the mappings we have used in 
  this article, for the exterior algebra, or the mappings of
  \cite{Sturmfels:LiftGB}, for the symmetric algebra case.

  Even in these cases, the generic initial ideal need not be
  Borel-fixed. 

  \begin{example}
    Suppose that \(V\) is 3-dimensional.
    Let \(I\) be the principal ideal (in the exterior algebra on
    \(V\)) on a generic quadratic form. Then 
    \(I\) is already in generic coordinates, and \(\init_{\prec_e}(I)
    = (x_2x_3)\), as long as \(x_1 \prec_e x_2 \prec_e x_3\).
    Hence, \(\gin(J)\) is generated by 
    \begin{displaymath}
      X_2X_3, \, X_1^2, \, X_2^2, \, X_3^2, \,
      X_2 X_1, \, X_3X_1, \, X_3X_2.
    \end{displaymath}
    Now let \(b\) the Borel transformation \(X_1 \mapsto X_1 +X_2\).
    We get that  \(b(X_1^2)=X_1^2+X_1X_2+X_2X_1+X_2^2\), which is not
    in the initial ideal.
  \end{example}

  \begin{example}
    Consider the commutator ideal in \(K\left \langle X_1, X_2 \right
    \rangle\). It is invariant under all linear coordinate changes, so 
    the generic initial ideal coincides with the initial ideal. The
    gin in degree \(2\) is  therefore \(X_2 X_1\), which is not Borel-fixed.
  \end{example}
\end{section}

\bibliographystyle{amsplain}
\bibliography{journals,articles}
\end{document}